\documentclass{article}
\usepackage[dvipsnames]{xcolor}
\usepackage{esvect}
\usepackage{graphicx}      
\usepackage{tikz}
\usepackage{amsmath,amssymb, amsfonts,amsthm}
\usepackage{subfig}
\usepackage{graphicx}
\newcommand{\R}{\mathbb{R}}

\newcommand{\cS}{\mathcal{S}}

\newcommand{\cX}{\mathcal{X}}

\newcommand{\supp}{\operatorname{supp}}

\newcommand{\rint}{\operatorname{int}}
\newcommand{\rL}{\mathrm{L}}

\usepackage{microtype}

\newcommand{\conv}{\operatorname{conv}}

\newcommand{\rrd}{\mathrm{d}}

\makeatletter
\def\blfootnote{\xdef\@thefnmark{}\@footnotetext}
\makeatother

\newcommand{\blue}[1]{\textcolor{blue}{#1}}

\newcommand{\xc}[1]{\vspace{.1cm}

\noindent {\em #1} }

\newcommand{\rLs}{\rL_s^\infty([0,1]^2,\R_{\geq 0})}

\newcommand{\rLsn}{\rL_s^\infty([0,1]^2,\R)}
\newcommand{\rLosn}{\rL^\infty([0,1],\R)}
\newcommand{\rLos}{\rL^\infty([0,1],\R_{\geq 0})}
\newtheorem{definition}{Definition}
\newtheorem{theorem}{Theorem}
\newtheorem{proposition}[theorem]{Proposition}

\usepackage[style= numeric-comp,hyperref=true, doi=false,url=false,
            isbn=false,
            firstinits=true, sorting = none, 
            block=none, backend=bibtex,maxnames=99]{biblatex}
\renewbibmacro{in:}{} 
\bibliography{bibliog}

\title{Graphons and the $H$-property}

\date{}
\author{Mohamed-Ali Belabbas\thanks{M.-A. Belabbas is with the Coordinated Science Laboratory, University of Illinois, Urbana-Champaign. Email: \texttt{belabbas@illinois.edu}} \quad and \quad Xudong Chen\thanks{X. Chen is with the Electrical and Systems Engineering, Washington University in St. Louis. Email: \texttt{cxudong@wustl.edu}. 
}
}

\begin{document}
\maketitle

\begin{abstract}               

\blfootnote{M.-A. Belabbas and X. Chen contributed equally to the manuscript in all categories.}
A graphon satisfies  the $H$-property if graphs sampled from it contain a Hamiltonian decomposition almost surely, which in turn implies that the corresponding network topologies are, e.g., structurally stable and structurally ensemble controllable. In recent papers~\cite{belabbas2021h,belabbas2023geometric}, we have exhibited a set of conditions that is essentially necessary and sufficient for the $H$-property to hold for the finite-dimensional class of step-graphons. The extension to the infinite-dimensional case of general graphons was hindered by the fact that said conditions relied on objects that do not admit immediate extensions to the infinite-dimensional case. We outline here our approach to bypass this difficulty and state conditions that guarantee that the $H$-property holds for general graphons.
\end{abstract}

%===============================================================================

\section{Introduction and Problem Formulation}

Structural system theory deals with the problem of understanding when a given network topology can sustain a prescribed  system property. Typical such properties are, e.g, controllability and stability. In more detail, consider a network of $n$ mobile agents $x_1,\ldots, x_n$, whose communication topology is described by a directed graph ({\em digraph}) $G = (V, E)$, with the nodes $v_1,\ldots, v_n$ representing the agents and directed edges $v_iv_j$ indicating the information flow (with the convention that a directed edge $v_iv_j$ indicates that agent $x_j$ can access state information from $x_i$). 
Given the digraph $G$, a system dynamics $\dot x(t) = f(x(t))$ is said to be {\em adapted to $G$} if the dynamics of $x_i(t)$ depend only on its incoming neighbors in $G$: $\frac{\partial f_i}{\partial x_j} \neq 0 \Rightarrow v_jv_i \in E$. 

We denote by $\Sigma_G$ the set of  differentiable dynamics adapted to $G$. 
Next, given a desired system property~$\cS$ (e.g., asymptotic stability at a given equilibrium), we say that the digraph $G$ {\em sustains} $\cS$ if there exists a dynamics $f \in \Sigma_G$ satisfying the property~$\cS$. This research line was initiated by C.-T. Lin in his seminal paper~\cite{lin1974structural} and led to several extensions, among which we mention the work of the authors on structural stability~\cite{belabbas2013sparse} and structural ensemble controllability~\cite{chen2021sparse} as they are the basis of the present work.

Graphons have recently been introduced in~\cite{lovasz2006limits, borgs2008convergent} to study large graphs. A graphon can be seen as both the limit object of a convergent  sequence (where convergence is in the cut-norm~\cite{frieze1999quick}) of graphs of increasing size, and as a statistical model from which to sample random graphs (we elaborate on this  below). 
Mathematically, a graphon is a symmetric, measurable function $W: [0,1]^2\to [0,1]$, 
with  $W(x,y)=W(y,x)$  for all $x,y\in [0,1]^2$. 
We note that a graphon should be viewed as an equivalence class of such functions, where $W_1 \sim W_2$ if there exists a measure-preserving map $\phi:[0,1] \to [0,1]$ such that $W_1(x,y) = W_2(\phi(x),\phi(y))$. This is the continuous equivalent of saying that graphs that are equal up to relabeling of their nodes are in the same equivalence class. We will overlook this distinction for the sake of clarity, but all our statements lift to the equivalence class.

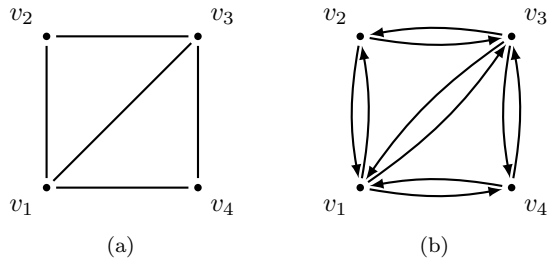
\begin{figure}
    \centering
    \subfloat[\label{sfig1:Ggraph}]{
    \begin{tikzpicture}[scale=1.15]
    \tikzset{every loop/.style={}}
		\node [circle,fill=black,inner sep=1pt,label=below left:{$v_1$}] (1) at (0, 0) {};
		\node [circle,fill=black,inner sep=1pt,label=above left:{$v_2$}] (2) at (0, 1.75) {};
		\node [circle,fill=black,inner sep=1pt,label=above right:{$v_3$}] (3) at (1.75, 1.75) {};
		\node [circle,fill=black,inner sep=1pt,label=below right:{$v_4$}] (4) at (1.75, 0) {};

  \path[draw,thick,shorten >=2pt,shorten <=2pt]
		 (1) edge[black] (2)
		 (1) edge[black] (4)
		 (2) edge[black] (3)
		 (3) edge[black] (4)
		 (3) edge[black] (1)
		 ;
\end{tikzpicture}
}
\qquad
\subfloat[\label{sfig1:vecG}]{
  \begin{tikzpicture}[scale=1.15]
    \tikzset{every loop/.style={}}
		\node [circle,fill=black,inner sep=1pt,label=below left:{$v_1$}] (1) at (0, 0) {};
		\node [circle,fill=black,inner sep=1pt,label=above left:{$v_2$}] (2) at (0, 1.75) {};
		\node [circle,fill=black,inner sep=1pt,label=above right:{$v_3$}] (3) at (1.75, 1.75) {};
		\node [circle,fill=black,inner sep=1pt,label=below right:{$v_4$}] (4) at (1.75, 0) {};
		
  \path[draw,thick,shorten >=2pt,shorten <=2pt]

		 (1) edge[,bend right=10,-latex] (2)
		 (2) edge[,bend right=10,-latex] (1)
 		 (1) edge[,bend right=10,-latex] (4)
		 (4) edge[,bend right=10,-latex] (1)
		 (3) edge[,bend right=10,-latex] (2)
		 (2) edge[,bend right=10,-latex] (3)
		 (3) edge[,bend right=10,-latex] (4)
		 (4) edge[,bend right=10,-latex] (3)

		 (3) edge[bend right=10,-latex] (1)
		 (1) edge[bend right=10,-latex] (3)
		 ;

\end{tikzpicture}
}
\caption{
{\em Left:} An undirected graph on $4$ nodes. {\em Right:} Its directed counterpart by replacing every undirected edge with two oppositely oriented edges.
}\label{fig:hamildecom}
\end{figure}

We next describe how to sample a graph $G_n$ on $n$ nodes from a graphon $W$. 
\xc{Sampling procedure} Let $\mathrm{Uni}[0,1]$ be the uniform distribution on $[0,1]$. Given a graphon $W$, a graph $G_n=(V,E)$ on  $n$ nodes sampled from~$W$, written as $G_n\sim W$, is obtained as follows: 
\begin{enumerate}
    \item Sample $x_1,\ldots,x_n\sim \mathrm{Uni}[0,1]$ independently. 
    We call $x_i$ the {\em coordinate of node} $v_i\in V$.\label{item:S1}

    \item For any two distinct nodes $v_i$ and $v_j$, place an edge $(v_i,v_j) \in E$ with probability $W(x_i,x_j)$.\label{item:S2}
\end{enumerate}
Note that $G_n$ is {\em undirected} and we denote by $\vec{G}_n=(V,\vec E)$ the {\em directed} version of $G_n$, with the edge set
$\vec E :=\{v_iv_j, v_jv_i \mid (v_i,v_j) \in E \}.$
In words, we replace an undirected edge $(v_i,v_j)$ with the pair of directed edges $v_iv_j$ and $v_jv_i$.  See Figure~\ref{fig:hamildecom} for illustration.

In this abstract, we characterize the graphons $W$ with the property that network topologies sampled from $W$ are structurally stable~\cite{belabbas2013sparse} or structurally ensemble controllable~\cite{chen2021sparse} with probability one. 
Precisely, we have the following definition~\cite{belabbas2023geometric}:

\begin{definition}[$H$-property]\label{def:Hproperty}
Let $W$ be a graphon and $G_n\sim W$. 
Then, $W$ has the {\bf $H$-property} if 
$$
\lim_{n\to\infty}\mathbb{P}(\vec G_n \mbox{ has a Hamiltonian decomposition}) = 1.
$$
\end{definition}
The problem we address can now be formulated as characterizing the $H$-property in graphons. 

\section{Results for Step-graphons}\label{sec:Stepg}
As a stepping-stone to the study of the general case of graphons, we have studied in our recent work~\cite{belabbas2023geometric} the class of {\em step-graphons} (see Definition~\ref{def:stepgraphon} below), and  provided essentially necessary and sufficient conditions  for the $H$-property to hold in that context. We summarize below the key objects  introduced to solve the problem as well as the main statement, as they are a blueprint for the study of the general case.

We start with the following definition: 

\begin{definition}[Step-graphon and its partition]\label{def:stepgraphon}
A graphon $W$ is a {\bf step-graphon} if there exists an increasing sequence $0 = \sigma_0 < \sigma_1< \cdots < \sigma_q = 1$, for $1 \leq q < \infty$, such that $W$ is constant over each rectangle $[\sigma_{i}, \sigma_{i + 1})\times [\sigma_{j}, \sigma_{j + 1})$ for all $0\leq i, j\leq q-1$.  The sequence $ \sigma = (\sigma_0,\sigma_1,\ldots,\sigma_q)$ is called a {\bf partition for $W$}.
\end{definition}

It should be clear that for a given step-graphon $W$, there exists infinitely many partition sequences (obtained, e.g., by sub-dividing intervals where the step-graphon is constant).

Step-graphons are related to stochastic block models~\cite{holland1983stochastic}, with the difference that nodes in stochastic block models are deterministically placed in intervals (called communities in that context) whereas, in the graphon case, the nodes are randomly placed.

Through our extant work, we have highlighted  three key objects associated with a step-graphon $W$ as being of particular importance.  
We introduce them below:

\begin{definition}[Concentration vector]\label{def:convec}
    For a step-graphon with partition $\sigma = (\sigma_0,\ldots,\sigma_q)$,  
    its {\bf concentration vector} is $x^* = (x^*_1,\ldots, x^*_q)$, where  
    $x^*_i := \sigma_i - \sigma_{i-1}$, for all $i = 1,\ldots, q$. 
\end{definition}

It should be clear from the sampling procedure given above that the concentration vector describes the {\em expected proportion} of sampled nodes in each interval. 
The support of a step-graphon can be described by a graph,  called {\em skeleton graph}, which is defined as

\sloppy
\begin{definition}[Skeleton graph]\label{def:skeleton}
To a step-graphon $W$ with partition $\sigma = (\sigma_0,\ldots, \sigma_q)$, we assign the  undirected graph $S = (U, F)$ on $q$ nodes, with $U =\{u_1,\ldots, u_q\}$ and edge set $F$ defined as follows: there is an edge between $u_i$ and $u_j$ if and only if $W$ is non-zero over $[\sigma_{i-1},\sigma_i)\times [\sigma_{j - 1}, \sigma_j)$. 
We call $S$ the {\bf skeleton graph} of $W$ for the partition sequence~$\sigma$.
\end{definition}
\fussy

\begin{figure}
\centering
\begin{tikzpicture}[scale=.4]
\filldraw [fill=Gray!70!black!40, draw=Gray!70!black!40] (1.2,0) rectangle (2.4,1.6);
\filldraw [fill=Gray!70!black!40, draw=Gray!70!black!40] (2.4,1.6) rectangle (4,2.8);
\filldraw [fill=black, draw=black] (2.4,0) rectangle (4,1.6);
\filldraw [fill=Gray!60!black!70, draw=Gray!60!black!70] (0,1.6) rectangle (1.2,2.8);
\filldraw [fill=Gray!60!black!70, draw=Gray!60!black!70] (1.2,2.8) rectangle (2.4,4);
\filldraw [fill=black, draw=black] (0,2.8) rectangle (1.2,4);

\draw [draw=black,very thick] (0,0) rectangle (4,4);
\node [] at (2, -1.2) {$W$}; \node [] at (5,2) {};
\end{tikzpicture}
\begin{tikzpicture}[scale=1.0]
    \tikzset{every loop/.style={}}
    \node [] (0) at (1,-1) {$S$};
		\node [circle,fill=black,inner sep=1.2pt,label=below:{$u_1$}] (1) at (0, 0) {};
		\node [circle,fill=black,inner sep=1.2pt,label=below:{$u_2$}] (2) at (.8, 0) {};
		\node [circle,fill=black,inner sep=1.2pt,label=below:{$u_3$}] (3) at (1.6, 0) {};
		
	\path[draw,thick,shorten >=2pt,shorten <=2pt]
		 (1) edge[black, loop left] (1)
          (3) edge[black, loop right] (3)
		 (1) edge[black] (2)
		 (2) edge[black] (3)
		 ;
\end{tikzpicture}
\begin{tikzpicture}[scale=1.1]
    \tikzset{every loop/.style={}}
    \node [circle,fill=black,inner sep=0.8pt] (0) at (0,0) {};
    \draw[->] (0,0) -- (1.55,0) node [above] {\footnotesize $x_1$};
    \draw[->] (0,0) -- (0,1.55) node [right] {\footnotesize $x_2$};
    \draw[->] (0,0) -- (-0.9,-0.9) node [left] {\footnotesize $x_3$};

    \node (1) at (1.2, 0) {};
    \node (2) at (0, 1.2) {};
    \node (3) at (-0.66, -0.66) {};
    \node (12) at (0.6, 0.6) {};
    \node (23) at (-0.33, 0.27) {};

\filldraw[draw=black, fill=gray!0]  (1.center) -- (2.center) -- (3.center) -- (1.center) -- cycle;

\filldraw[draw=blue!20, fill=blue!20]  (1.center) -- (12.center) -- (23.center) -- (3.center) -- cycle;

\node [] (xs) at (0.8,1.0) {\small \blue{$\cX(S)$}};
\node [] (xsl) at (0.8,1.0) {};
\node [] (xst) at (0.7, 0.1) {};

\path[draw,shorten >=3pt,shorten <= 3pt]
    (xsl) edge[-latex] (xst);

\node [circle,fill=red,inner sep=1.0pt,label=below:{\color{red} $x^*$}] (x) at (0.096,0.096) {};
\end{tikzpicture}
\caption{{\em Left:} A step-graphon $W$ with partition $\sigma = (0,0.3,0.6,1)$. {\em Middle:} Skeleton graph $S$. {\em Right:} Edge polytope $\cX(S)$ and the concentration vector $x^*=(0.3,0.3,0.4)$.}\label{fig:edgepoly}
\end{figure}
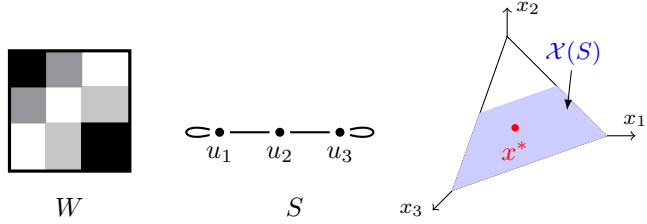

The last object we introduced is derived from the skeleton graph.  
Let $F=\{f_1,\ldots, f_r\}$, we let the edge-incidence matrix  $B \in \R^{q \times r}$ of a graph $S=(U,F)$ be given as
\begin{equation}\label{eq:defZS}
b_{ij} := \frac{1}{2}
\begin{cases}
    2, & \text{if } f_j\in F \text{ is a self-loop on node } u_i,       \\
    1, & \text{if node } u_i \text{ is incident to } f_j\in F, \\
    0, & \text{otherwise}.
\end{cases}
\end{equation}
Owing to the factor $\frac{1}{2}$ in~\eqref{eq:defZS}, all columns of $B$ are probability vectors, i.e., all entries are nonnegative and sum to one. The edge polytope of $S$, introduced in~\cite{ohsugi1998normal}, is 
\begin{definition}[Edge polytope]\label{def:edgepolytope}
Let $S = (U,F)$ be a skeleton graph and $B$ be the associated incidence matrix. Let $b_j$, for $1\leq j \leq |F|$, be the columns of $B$.  
The {\em edge polytope} of $S$, denoted by $\cX(S)$, is the convex hull generated by the $b_j$'s: $\cX(S):= \conv\{b_j \mid j = 1,\ldots, r\}$.

\end{definition}
We illustrate the above three key objects in Figure~\ref{fig:edgepoly}.  

To state the result, we introduce the following three conditions:

\begin{description}
\item[\it Condition $A$:] The skeleton graph $S$ has an odd cycle.
\item[\it Condition $B$:] The concentration vector $x^*$ belongs to the {\em relative interior} of the edge polytope $\cX(S)$, i.e.,  $x^* \in \rint \cX(S)$.
\item[\it Condition $B'$:] $x^* \in \cX(S)$.
\end{description}

The above conditions can be shown to be independent of the choice of a partition sequence~$\sigma$~\cite{belabbas2023geometric}. 
The main results proved in~\cite{belabbas2021h,belabbas2023geometric} are reproduced below
\begin{theorem}\label{thm:Hproperty}
For almost all step-graphons $W$, the probability that $\vec G_n\sim W$ has a Hamiltonian decomposition tends to either $0$ or $1$. Furthermore, 
\begin{itemize}
\item If Conditions $A$ and $B$ hold, then $W$ has the $H$-property.
\item If either Condition $A$ or $B'$ does not hold, then $W$ does not have the $H$-property. In fact, for this case, 
$$
\lim_{n\to\infty}\mathbb{P}(\vec G_n \mbox{ has a Hamiltonian decomposition}) = 0.
$$
\end{itemize}
\end{theorem}
For example, the step-graphon in Figure~\ref{fig:edgepoly} satisfies Conditions A and B and, hence, has the $H$-property.

\section{Extension to General Graphons}
We now discuss the general case of graphons. One first observes that while some objects introduced above, such as the concentration vector, admit a relatively simple translation to the infinite-dimensional case, the skeleton graph does not admit a simple infinite-dimensional equivalent. Since the edge-polytope was defined from the edge-incidence matrix of the skeleton graph, naive limits $n \to \infty$ do not yield useful results.

To circumvent this difficulty, we are forced to take a different approach to describing the edge-polytope of a graphon, that by-passes the use of the skeleton graph. Relying on some technical results in~\cite{belabbas2023geometric}---precisely, what we refer to as the $A$-matrix of a step-graphon---we arrive at the following characterizations:

\xc{Extension of Condition A.} 
Let $\rLsn$ be the space of all measurable, bounded, and symmetric functions from $[0,1]^2$ to $\R$.  
To a given graphon $W$, we introduce the map:
\begin{equation}\label{eq:defphiw}
\Phi_{W}: \rLsn \to \rLosn:  c \mapsto x(s):=\int^1_0 W(s,t) c(s,t) \rrd t.
\end{equation}
Also, we let $\overline W$ be the {\em saturation} of $W$, which is a graphon valued in $\{0,1\}$ defined as
$$
\overline W(x,y) := 
\begin{cases}
1 & \mbox{if } W(x,y) > 0, \\
0 & \mbox{otherwise}.
\end{cases}    
$$
Graphons valued in $\{0,1\}$ are also called {\em random-free graphons}~\cite{janson2010graphons}.

The following result, which characterizes the map $\Phi(W)$ when $W$ is a step-graphon,  highlights its  importance to the analysis of the $H$-property.

\begin{proposition}\label{prop:surjectivity}
For $W$ a step-graphon and $S$ its skeleton graph, the map $\Phi_{\overline W}$ is surjective {\em if and only if} $S$ has an odd-cycle.
\end{proposition}

We thus introduce the following condition as an extension of Condition~$A$:

\begin{description}
\item[\it Condition $A_{\rm ext}$:] $\Phi_{\overline W}$ is surjective. 
\end{description}

\xc{Extension of Condition B.} To proceed, we introduce the operator that integrates a function $x(s)$ over each interval of a given partition $\sigma =(\sigma_0,\ldots,\sigma_q)$:  
\begin{equation}\label{eq:defmu}
\mu_\sigma:\rLosn \to \R^q: x \mapsto \bar x:=\begin{bmatrix}\displaystyle\int_{\sigma_0}^{\sigma_1} x(s) \rrd s \\ \vdots \\ \displaystyle \int_{\sigma_{q-1}}^{\sigma_q} x(s) \rrd s
\end{bmatrix}.
\end{equation}

For $W$ a step-graphon with skeleton graph $S$, we define
$$
\rL(W):=\{x(s) \in \rLos \mid \mu_\sigma(x)  \in  \cX(S)\}.
$$
We state, without a proof, that $\rL(W)$ does not depend on the choice of partition sequence $\sigma$ (and hence, of $S$) for $W$.  

The space $\rL(W)$ can be thought of as a functional equivalent of the edge-polytope of a step-graphon, as the next proposition shows:
\begin{proposition}\label{prop:prorl}
     Let $\mathbf{1}\in \rLos$ be the unit constant function. For $W$  a step-graphon with concentration vector $x^*$ and skeleton graph $S$,  $\mathbf{1}\in \rL(W)$ (resp. $\mathbf{1}\in \rint \rL(W)$) {\em if and only if} $x^*\in \cX(S)$ (resp. $x^*\in \rint \cX(S)$).
\end{proposition}

Even though $\rL(W)$ is a functional space, it still relies on $S$ for its definition and thus cannot be applied directly to general $W$. 
To resolve the issue, we define 
\begin{equation}\label{eq:defcxw}
\cX(W):= \left \{\Phi_{\overline W}(c)  \bigm|  c \in \rLs  \mbox{ and } 
\|\Phi_{\overline W}(c)\|_{\mathrm{L}^1} = 1\right \}.
\end{equation}
Note that the requirement $\|\Phi_{\overline W}(c)\|_{\mathrm{L}^1} = 1$ is equivalent to $ \int_{\supp W} c(s,t) \rrd s\rrd t = 1$. 
It should be clear that $\cX(W)$ is well defined for all graphons. 

We have the following result:

\begin{proposition}\label{prop:rlwxs}
For $W$ a step-graphon, $\cX(W)=\rL(W)$. 
\end{proposition}

This proposition states that the $\cX(W)$ introduced above, which did not rely on the skeleton graph of $W$, equals the space $\rL(W)$ and together with Proposition~\ref{prop:prorl} shows that it can indeed serve as an equivalent of the edge-polytope applicable to the case of general graphons. This leads to the following extension of Condition B:

\begin{description}
\item[\it Condition $B_{\rm ext}$:]  $\mathbf{1}\in \rint \cX(W)$. 
\end{description}

We are now in a position to state the main result of this abstract: 

\begin{theorem}\label{thm:hgraphon}
Given a graphon $W$, let $\Phi_{\overline W}$ and $\cX(W)$ be as  in Eqns.~\eqref{eq:defphiw} and~\eqref{eq:defcxw}, respectively. Then, $W$ has the $H$-property if Conditions $A_{\rm ext}$ and $B_{\rm ext}$ are satisfied.
\end{theorem}

The proofs of the above results will appear in an upcoming paper. 

\printbibliography
   
\end{document}